\documentclass{amsart}
\usepackage{amssymb, latexsym}

 \newtheorem{thm}{Theorem}[section]
 
 \newtheorem{cor}[thm]{Corollary}
 \newtheorem{lem}[thm]{Lemma}
 \newtheorem{prop}[thm]{Proposition}
 
 \theoremstyle{definition}
 \newtheorem{defn}[thm]{Definition}
 \theoremstyle{remark}
 \newtheorem{rem}[thm]{Remark}
 \theoremstyle{example}
 \newtheorem{ex}[thm]{Example}
 \numberwithin{equation}{section}

 \newcommand{\To}{\rightarrow}
 \newcommand{\bP}{\mathbb{P}}
 \newcommand{\bC}{\mathbb{C}}
 \newcommand{\bZ}{\mathbb{Z}}
 
 \newcommand{\bQ}{\mathbb{Q}}
 
 \newcommand{\DM}{Deligne-Mumford}
 \newcommand{\GW}{Gromov-Witten}
 \newcommand{\sX}{\mathcal{X}}
 \newcommand{\sC}{\mathcal{C}}
 \newcommand{\sK}{\mathcal{K}}
 \newcommand{\sM}{\mathcal{M}}

 \newcommand{\sO}{\mathcal{O}}
 \newcommand{\ii}{\mathcal{I}}
 \newcommand{\ri}{\bar{\mathcal{I}}}

 \newcommand{\generalK}{\mathcal{K}_{g,n}(\sX,\beta)}
 \newcommand{\generalKp}{\mathcal{K}_{g,n+1}(\sX,\beta)}
 \newcommand{\coloneq}{\mathrel{\mathop:}=}
 \newcommand{\e}{e_{1, \cdots, n}}
 \newcommand{\ep}{e_{1, \cdots, n+1}}
 \newcommand{\f}{\check{e}_{n+1}}
 \newcommand{\TABLESIZE}{.92}

\usepackage{stmaryrd}
\usepackage{diagrams}
\newarrow{Eq}=====
\newarrow{Into}C--->
\newarrow{Mapsto}|--->

\begin{document}

\title[A Reconstruction Theorem for Genus Zero Gromov-Witten Invariants of Stacks]
 {A Reconstruction Theorem for Genus Zero Gromov-Witten Invariants of Stacks}

\author{ Michael A. Rose }

\address{Department of Mathematics, University of Wisconsin-Madison, Madison, Wisconsin, 53706,
 United States}

\email{rose@math.wisc.edu}


\begin{abstract}
We generalize the \emph{First Reconstruction Theorem} of Kontsevich and Manin in two respects.  First, we allow the target space to be a \DM\ stack.  Second, under some convergence assumptions, we show it suffices to check the hypothesis of $H^2$-generation not on the cohomology ring, but on an any quantum ring in the family given by small quantum cohomology.  As an example the latter result is used to compute genus zero \GW\ invariants of $\bP(1,b)$.
\end{abstract}

\maketitle

\section*{Introduction}

The main goal of this article is to prove an appropriate extension of the \emph{First Reconstruction Theorem} of Kontsevich and Manin \cite[Theorem 3.1]{KM} in genus zero \GW\ theory to the case where the target is a \DM\ stack.  More precisely, consider the following version of [ibid]:

\begin{thm} \label{KM_reconstruction}
Let $X$ be a smooth projective variety.  Suppose that $H^*(X)$ is generated by $H^2(X)$, then all genus zero \GW\ invariants can be uniquely reconstructed starting from 3-point invariants.
\end{thm}

The proof relies mainly on the \emph{Witten-Dijkgraaf-Verlinde-Verlinde} (WDVV) \cite{Witten,DVV1,DVV2} equations in genus zero \GW\ theory and uses the hypothesis in two ways.  First, note that the divisor axiom applies precisely to classes in $H^2(X)$.  Second, note that the degree zero 3-point invariants on $X$ form the structure constants of the cup product on $H^*(X)$, so the hypothesis can be viewed as a condition on degree zero 3-point invariants.

Now, let $\sX$ be a smooth \DM\ stack with projective coarse moduli space $X$.  \GW\ invariants on $\sX$ pull back classes from $H^*(\ri_{\mu}(\sX))$, where $\ri_{\mu}(\sX)$ is a stack naturally associated to $\sX$ called the \emph{rigidified cyclotomic inertia stack}.  $\ri_{\mu}(\sX)$ contains $\sX$ as an open and closed substack, and the subspace $H^*(\sX) \subseteq H^*(\ri_{\mu}(\sX))$ is called the \emph{untwisted sector}.  To see how Theorem \ref{KM_reconstruction} might be generalized to this context consider the following results.  In \cite{A} and \cite{AGV2}, the WDVV equation is extended to the case of \DM\ stacks.  There is a divisor axiom which in this context applies to classes in $H^2(\sX) \subset H^*(\ri_{\mu}(\sX))$, and degree zero 3-point invariants determine a new product on the vector space $H^*(\ri_{\mu}(\sX))$.  The new ring is called the \emph{orbifold cohomology ring} of $\sX$ and is denoted $H^*_{orb}(\sX)$ \cite{CR}.  The following proposition then follows via the same technique as the proof of Theorem \ref{KM_reconstruction}.
\begin{prop}
Let $\sX$ be a smooth \DM\ stack with projective coarse moduli space.  Suppose that $H^*_{orb}(\sX)$ is generated by $H^2(\sX)$, then all genus zero \GW\ invariants can be reconstructed starting from 3-point invariants.
\end{prop}

However, this is only a theoretical generalization:  unless $\sX$ is a scheme $H^*_{orb}(\sX)$ is never generated by $H^2(\sX)$.  Thus we search for a more useful generalization.  When the stack $\sX$ satisfies $\text{deg } c_1(T_{\sX}) > 0$ (or more generally satisfies a condition guaranteeing the convergence of the small quantum product, see Definition \ref{conv_crit}), $H^*_{orb}(\sX)$ is a specialization in a family of rings given by the small quantum ring $QH^*(\sX)$.  All the rings have the same underlying vector space, namely $H^*(\ri_{\mu}(\sX))$, and we may weaken the above hypothesis by requiring $H^2(\sX)$-generation not on $H^*_{orb}(\sX)$, but on any one of the specializations in the family.  If $q_1, \cdots, q_d$ denote parameters on the base of the family and $\bar{\lambda} = (\lambda_1, \cdots, \lambda_d)$, let $QH^*_{\bar{\lambda}}(\sX)$ denote the specialization of $QH^*(\sX)$ given by setting $q_i = \lambda_i$.  The following is then the main theorem of this article.

\begin{thm} \label{intro_reconstruction}
Let $\sX$ be a smooth \DM\ stack with projective coarse moduli space satisfying the \emph{convergence criterion} (Definition \ref{conv_crit}).  Suppose for some $\bar{\lambda}$, $QH_{\bar{\lambda}}(\sX)$ is generated by $H^2(\sX)$, then all genus zero \GW\ invariants can be reconstructed starting from 3-point invariants.
\end{thm}

The paper is organized as follows.  In Section 1, we recall the basic theory of \GW\ invariants of \DM\ stacks as developed in \cite{AGV} and \cite{AGV2}.  Section 2 contains the conventions and definitions of small quantum cohomology and its specializations needed in the proof of Theorem \ref{intro_reconstruction} which is given in Section 3.  Finally, Section 4 contains an application of Theorem \ref{intro_reconstruction} to the stacks $\bP(1,b)$.  Here the proof gives an explicit algorithm for computing arbitrary genus zero invariants of $\bP(1,b)$.  Example computations for low values of $b$ are tabulated in the appendix.

Many thanks are due to my advisor Lev Borisov for many stimulating conversations.  I would also like to thank Dan Abramovich for information on the literature of \GW\ theory of stacks, and Charles Cadman and Hiroshi Irritani for their useful comments.


\section{Gromov-Witten Theory of \DM\ Stacks}

Assume all stacks and schemes lie over the complex numbers.

The \GW\ theory of \DM\ stacks was first mathematically formulated by Chen and Ruan \cite{CR} in the symplectic setting and subsequently by Abramovich, Graber, and Vistoli \cite{AGV,AGV2} in the algebraic setting.  In this article, we follow the latter approach.  See also \cite{A} and \cite{AV1}.

When attempting to construct some analog of a Kontsevich moduli space but with target a \DM\ stack, a natural way to get a compact space is to allow the domain curves to acquire stack structure.  The development of the theory then begins by carefully characterizing these curves.  For simplicity, we give the rigorous definition only over a point.

\subsection{Twisted curves and twisted stable maps}

\begin{defn}
An \emph{n-pointed twisted nodal curve} $(\sC, \Sigma_1, \cdots, \Sigma_n)$ is a diagram
\begin{diagram}
\sqcup_{i=1}^n \Sigma_i & \rInto & \sC \\
                        &        & \dTo_{\pi} \\
                        &        & C \\
\end{diagram}
where
\begin{enumerate}
  \item $\sC$ is a proper \DM\ stack with coarse moduli scheme $C$
  \item $(C, \pi(\Sigma_1), \cdots, \pi(\Sigma))$ is a n-pointed nodal curve
  \item Over the node of $C$ with local expression $\{xy = 0\}$, $\sC$ has \'{e}tale chart $$[\{xy = 0\} / \mu_r]$$ where the action is given by $(x,y) \mapsto (\xi u, \xi^{-1} v)$
  \item Over a marked point $\pi(\Sigma_i)$ of $C$, $\sC$ has \'{e}tale chart $$[\mathbb{A}^1 / \mu_r]$$ where the action is given by $u \mapsto \xi u$ and $\Sigma_i$ is the substack defined by $u = 0.$
\end{enumerate}

After appropriately defining n-pointed twisted nodal curves (and morphisms of n-pointed twisted nodal curves) over an arbitrary base scheme, we have the following theorem.

\begin{thm} [\cite{O}, Theorem 1.9]
The category $\bar{\mathcal{M}}_{g,n}^{tw}$ of n-pointed twisted nodal curves is a smooth Artin stack of dimension 3g - 3 + n.
\end{thm}

\end{defn}
Let $\sX$ be a \DM\ stack with projective coarse moduli scheme $X$.
\begin{defn}
An \emph{n-pointed twisted stable map} $(\sC \xrightarrow{f} \sX, \Sigma_1, \cdots, \Sigma_n)$ is a diagram
\begin{diagram}
\sqcup_{i=1}^n \Sigma_i & \rInto & \sC  & \rTo^{f}            & \sX          \\
                        &        & \dTo &                     &  \dTo        \\
                        &        & C    & \rTo^{\bar{f}} & X            \\
\end{diagram}
where
\begin{enumerate}
  \item $(\sC, \Sigma_1, \cdots, \Sigma_n)$ is an n-pointed twisted nodal curve
  \item $f$ is representable with $\bar{f}$ the induced map on coarse moduli spaces
  \item $\bar{f}$ is stable in the sense of Kontsevich \cite{Ko}.
\end{enumerate}
\end{defn}

Let $N_1(X)$ be the group of numerical equivalence classes of curves in $X$, and let $N^+(X) \coloneq N_1^+(X)$ be the monoid of effective classes.  Then for $\beta \in N_1(X)$ and for $g$ a non-negative integer, one says that $(\sC \xrightarrow{f} \sX, \Sigma_1, \cdots, \Sigma_n)$ has degree $\beta$ and genus $g$ if the stable map $\bar{f}$ does.

After appropriately defining n-pointed twisted stable maps (and morphisms of n-pointed twisted stable maps) over an arbitrary base scheme, we have the following theorem.

\begin{thm} [\cite{AV1}, Theorem 1.4.1] \label{K_proper}
    The category $\sK_{g,n}(\sX, \beta)$ of n-pointed twisted stable maps of genus g and degree $\beta$ is a proper \DM\ stack.  The coarse moduli space $\mathbf{K}_{g,n}(\sX, \beta)$ of $\sK_{g,n}(\sX, \beta)$ is projective.
\end{thm}


\subsection{Evaluation Maps and Inertia Stacks}

To evaluate an n-pointed twisted stable map $(\sC \xrightarrow{f} \sX, \Sigma_1, \cdots, \Sigma_n)$ at the $i^{th}$ substack $\Sigma_i$, one merely wants to restrict $f$ to obtain the representable morphism
\begin{diagram}
\Sigma_i  & \rTo^{f|_{\Sigma_i}} & \sX. \\
\end{diagram}
This morphism, however, contains more data than merely a point in $\sX$, so to keep track of this structure the evaluation should take values in the \emph{cyclotomic inertia stack}
\begin{alignat}{1}
\ii_{\mu}(\sX) &\coloneq \bigcup_r \ii_{\mu_r}(\sX) \notag \\
               &\coloneq \bigcup_r \text{\underline{Hom}}^{\text{rep}}(B(\mu_r), \sX) \label{homrep},
\end{alignat}
where the superscript \emph{rep} denotes the substack of representable morphisms (see \cite{O-Hom}).  Note that each $\Sigma_i \cong B(\mu_r)$ for some r.  However, there is difficulty in defining the evaluation on arbitrary families of twisted stable maps.  To remedy this we need to replace (\ref{homrep}) with the \emph{rigidified cyclotomic inertia stack}
\begin{equation} \label{gerbes}
\ri_{\mu}(\sX) \coloneq \bigcup_r \ri_{\mu_r}(\sX)
\end{equation}
where $\ri_{\mu_r}(\sX)$ is the stack classifying $\mu_r$-banded gerbes in $\sX$.  For base scheme $T$, an object of $\ri_{\mu_r}(\sX)(T)$ is a diagram
\begin{diagram}
\mathcal{G} & \rTo{f}   & \sX \\
\dTo_{p}    &           & \\
T           &           & \\
\end{diagram}
where $(\mathcal{G} \xrightarrow{p} T)$ is a gerbe banded by $\mu_r$ and $f$ is a representable morphism.

Since (\ref{homrep}) and (\ref{gerbes}) are compared to each other and to the inertia stack $\ii(\sX)$ in detail in \cite[Section~3]{AGV2}, we content ourselves with a couple of examples to illustrate their differences.

\begin{ex}
Let $Y$ be a smooth projective scheme, $G$ a finite group scheme, and let $\sX = [Y/G]$.  Then $$\ii_{\mu}(\sX) = \bigsqcup_{(g)} [Y^g/C(g)]$$ where the union is over conjugacy classes $(g)$ of $G$, $Y^g$ denotes the fixed locus of $g$, and $C(g)$ denotes the centralizer.  Note that the cyclic group $\langle g \rangle$ acts trivially on $Y^g$ and thus $C(g)/\langle g \rangle$ also acts on $Y^g$.  We then have $$\ri_{\mu}(\sX) = \bigsqcup_{(g)} [Y^g/(C(g)/\langle g \rangle)].$$
\end{ex}
\begin{ex} \label{mainexample}
Let $b \in \bZ_{>0}$, and let $\sX = \bP(1,b) \coloneq [(\mathbb{A}^2 \setminus \mathbf{0})\ /\ \bC^*]$ where the action is given by $(x,y) \mapsto (\xi x,\ \xi^b y)$.  We identify the only non-trivial isotropy group with the subgroup $\mu_b \subset \bC^*$.  Since $\bC^*$ is abelian we have  $$\ii_{\mu}(\sX) = \sX \sqcup \bigsqcup_{1 \neq g \in \mu_b} B(\mu_b).$$  Moreover, we have $$\ri_{\mu}(\sX) = \sX \sqcup \bigsqcup_{1 \neq g \in \mu_b} B(\mu_b / \langle g \rangle).$$
\end{ex}

Finally, we have evaluation morphisms
\begin{diagram}
\generalK & \rTo{e_i} & \ri_{\mu}(\sX) \\
(\sC \xrightarrow{f} \sX, \Sigma_1, \cdots, \Sigma_n) & \rMapsto & (\Sigma_i \xrightarrow{f|_{\Sigma_i}} \sX) \\
\end{diagram}
and the twisted evaluation morphisms defined by $\check{e}_i \coloneq i \circ e_i$ where $i$ is the isomorphism
\begin{diagram}
\ri_{\mu}(\sX) & \rTo{i} & \ri_{\mu}(\sX) \\
(\mathcal{G} \To \sX_S) & \rMapsto & ({}^{\tau}\mathcal{G} \To \sX_S). \\
\end{diagram}
Here $\tau: \mu_r \To \mu_r$ is the involution sending $\xi \mapsto \xi^{-1}$ and ${}^{\tau}\mathcal{G}$ denotes the gerbe induced by changing the $\mu_r$-banding on $\mathcal{G}$ by $\tau$.

In general $\ri_{\mu}(\sX)$ has many components.  If $\ri_{\mu}(\sX) = \bigsqcup_{i \in I} \sX_i$ is the decomposition into connected components, the evaluation morphisms induce a decomposition of $\generalK$ into open and closed substacks.  Writing
$$\sK_{g,n}(\sX, \beta, i_1, \cdots, i_n) \coloneq e_1^{-1}(\sX_{i_1}) \cap \cdots \cap e_n^{-1}(\sX_{i_n})$$
we then have
\begin{equation} \label{decomp}
\generalK = \bigsqcup_{(i_1, \cdots, i_n) \in I^n} \sK_{g,n}(\sX, \beta, i_1, \cdots, i_n).
\end{equation}


\subsection{Virtual fundamental classes}
\GW\ invariants require an intersection theory on $\generalK$.  While this stack is in general singular and not equidimensional, $\sK_{g,n}(\sX, \beta, i_1, \cdots, i_n)$ at least has a constant expected dimension.

Consider the following universal diagram where $\mathcal{U}$ is the universal twisted stable curve, and $f$ is the universal representable morphism:
\begin{diagram}
\mathcal{U} & \rTo^{f}  & \sX \\
\dTo^{\pi}  &           & \\
\sK_{g,n}(\sX, \beta, i_1, \cdots, i_n)   &       & \\
\end{diagram}

\begin{defn}
The \emph{expected dimension} of $\sK_{g,n}(\sX, \beta, i_1, \cdots, i_n)$, denoted $edim$ is the integer
$$edim = \chi(\sC, f^*T_{\sX}) + dim (\bar{\sM}_{g,n}^{tw})$$
where $\sC$ is any closed curve in the family.
\end{defn}

We compute the expected dimension using a Riemann-Roch theorem on twisted nodal curves.  Since the data of a vector bundle on a twisted curve also encodes the action of an isotropy group on a fiber, this action may also play a role in a Riemann-Roch formula.  This action is encoded via the notion of \emph{age}.

A representation $\rho: \mu_r \To \bC^*$ is determined by an integer $0 \le k \le r-1$ with $\rho(g) = g^k$.  Define
$$age(\rho) \coloneq \frac{k}{r}$$
This extends by linearity to a function on the representation ring
$$age: R\mu_r \To \bQ.$$
For any diagram in $\ri_{\mu_r}(\sX)$:
\begin{diagram}
\mathcal{G}     & \rTo^{f}  & \sX \\
\dTo            &           & \\
T               &           & \\
\end{diagram}
the gerbe $\mathcal{G} \To T$ is locally trivial.  The restriction of the vector bundle $f^*T_{\sX}$ to any trivialization defines a representation of $\mu_r$.  The age of this representation is independent of the trivialization, and we thus obtain a locally constant function also denoted by $age$:
$$age: \ri_{\mu}(\sX) \To \bQ.$$
The various uses of this notation will be clear from the context.  The Riemann-Roch theorem is then stated as follows:
\begin{thm} \cite{AGV2}
Let $\mathcal{E}$ be a locally-free sheaf on a twisted nodal curve $\mathcal{C}$.  For each point $p \in \mathcal{C}$, $\mathcal{E}|_p$ has a $G_p$-action where $G_p$ is the isotropy group at $p$ and $G_p \cong \mu_r$ for some $r$.  Then
$$\chi(\mathcal{E}) = rk(\mathcal{E})\chi(\sO_{\sC}) + \int_{\mathcal{C}} c_1(\mathcal{E}) - \sum_{p \in \sC} age(\mathcal{E}|_p).$$
\end{thm}

\begin{cor}
The expected dimension of $\sK_{g,n}(\sX, \beta, i_1, \cdots, i_n)$ is given by
$$edim = \int_{\beta} c_1(T_{\sX}) + dim(\sX)(1 - g) + 3g - 3 + n - \sum_{k=1}^n age(\sX_{i_k}).$$
\end{cor}

The techniques of \cite{Be} and \cite{BF} can be used to construct a virtual fundamental class
$$[\sK_{g,n}(\sX, \beta, i_1, \cdots, i_n)]^{vir} \in A_{edim}(\sK_{g,n}(\sX, \beta, i_1, \cdots, i_n))_{\bQ}.$$
These classes then give a virtual fundamental class on $\generalK$ by passing to the decomposition (\ref{decomp}):
$$[\generalK]^{vir} \in A_*(\generalK)_{\bQ}.$$


\subsection{Gromov-Witten classes}
We now define \GW\ classes.  There are several conventions that we may proceed with, but for simplicity we work with cohomology classes and we assume $X$ has only even cohomology.  Moreover we shall work with coefficients in some field $R$.

Let $\e \coloneq e_1 \times \dots \times e_n$ denote the product of the evaluation maps.  Consider the following diagram:
\begin{diagram}
\generalKp & \rTo^{\e} & (\ri_{\mu}(\sX))^n \\
\dTo_{\f}  &           &   \\
\ri_{\mu}(\sX)     &           &   \\
\end{diagram}
and let $r: \ri_{\mu}(\sX) \To \bZ$ be the locally constant function taking value $r$ on $\ri_{\mu_r}(\sX)$.  Denote also by $r$ the induced class in $H^0(\ri_{\mu}(\sX))$.
\begin{defn}
For classes $\gamma_1, \cdots , \gamma_{n} \in H^*(\ri_{\mu}(\sX))$, and $\beta \in N^+_1(X)$,
$$\langle \gamma_1, \cdots, \gamma_n, \ast \rangle_{g, \beta}^{\sX} \coloneq r \cdot (\f)_{*}(\e^*(\gamma_1 \times \dots \times \gamma_n) \cap [\mathcal{K}_{g,n}(\sX,\beta)]^{vir})$$ is a \emph{\GW\ cohomology class}.
\end{defn}
Moreover, define the following \emph{cohomological \GW\ numbers} in $R$:
\begin{defn}
For $\gamma_1, \cdots , \gamma_{n+1} \in H^*(\ri_{\mu}(\sX))$,
$$ \langle \gamma_1, \cdots, \gamma_n, \gamma_{n+1} \rangle_{g, \beta}^{\sX} \coloneq \int_{[\mathcal{K}_{g,n}(\sX,\beta)]^{vir}} \ep^*(\gamma_1 \times \dots \times \gamma_{n+1}).$$
\end{defn}

\newcommand{\g}{g_{ij}}
\newcommand{\gi}{g^{ij}}
These will frequently also be called \emph{\GW\ invariants} or \emph{(n+1)-point invariants}.
It will be convenient to interchange between the two formalisms in the following sections.  We see that the \GW\ classes and numbers determine each other by the following lemma whose proof is left to the reader.
\begin{lem} \label{lemma}
\begin{enumerate}
  \item{$$ \langle \gamma_1, \cdots, \gamma_n, \gamma_{n+1} \rangle_{g, \beta}^{\sX} = \int_{\ri_{\mu}(\sX)} r^{-1}\langle \gamma_1, \cdots, \gamma_n, \ast \rangle_{g, \beta}^{\sX}\ \cup\ i^*(\gamma_{n+1})$$}\label{1}

  \item{Let $\alpha_0, \cdots, \alpha_m$ be a basis for $H^*(\ri_{\mu}(\sX))$, let $\g = \int_{\ri_{\mu}(\sX)} r^{-1}\alpha_i \cup i^*(\alpha_j)$, let $(\gi)$ be the inverse matrix of $(\g)$, then $$\langle \gamma_1, \cdots, \gamma_n, \ast \rangle_{g, \beta}^{\sX} = \sum_{i,j}\langle \gamma_1, \cdots, \gamma_n, \alpha_i \rangle_{g, \beta}^{\sX}\ \gi\ \alpha_j.$$ }\label{2}
\end{enumerate}
\end{lem}

\subsection{Properties}

We introduce a couple more notational conventions to describe those properties of \GW\ invariants we will need.
\begin{defn}
The \emph{orbifold degree} of $\gamma \in H^*(\ri_{\mu_r}(\sX))$ is given by
$$orbdeg(\gamma) \coloneq deg(\gamma) + 2 age(\ri_{\mu_r}(\sX)).$$
\end{defn}
This defines a grading on the vector space $H^*(\ri_{\mu}(\sX))$, called the \emph{orbifold grading}.

Next, let $\sX \xrightarrow{p} X$ denote the morphism to coarse moduli space $X$.  For any $\beta \in N^+(X)$ and any $D \in H^2(\sX)$ we write
$$\int_{\beta} D\ \coloneq \int_{\beta'} D$$
where $\beta'$ is any class with $p_*(\beta') = \beta$.

When $\sX$ is a smooth projective variety, the following properties are well known and follow from \cite{Be}.  For $\sX$ a smooth \DM\ stack as above, the properties are proven in \cite[Sections~7,8]{AGV2}.

\begin{thm}
  \begin{enumerate}
    \item (Degree Axiom)
    For classes $\gamma_1, \cdots, \gamma_n \in H^*(\ri_{\mu}(\sX))$ homogeneous in the orbifold grading,
    $$\langle \gamma_1, \cdots, \gamma_n \rangle_{g, \beta}^{\sX} = 0$$
    unless
    $$\sum_i orbdeg(\gamma_i) = 2 \int_{\beta} c_1(T_{\sX}) \ + 2\ dim \sX (1 - g) + 2(3g - 3 + n).$$ \\
    \item (Fundamental Class Axiom)
    If $1_{\sX} \in H^*(\ri_{\mu}(\sX))$ denotes the fundamental class of the untwisted component $\sX \subseteq \ri_{\mu}(\sX)$, then $$\langle 1_{\sX}, \gamma_1, \cdots, \gamma_n \rangle_{g, \beta}^{\sX} = 0$$ unless $g = 0$ and $n = 2.$\\

    \item(Untwisted Divisor Axiom)
    For any $D \in H^2(\sX) \subset H^*(\ri_{\mu}(\sX))$, if $(\beta, g, n)$ is not any of $\beta = 0, g = 0, n < 3$ or $\beta = 0, g = 1, n =0$, $$\langle D, \gamma_1, \cdots, \gamma_n \rangle_{g, \beta}^{\sX} = (\int_{\beta} D) \cdot \langle \gamma_1, \cdots, \gamma_n \rangle_{g, \beta}^{\sX}.$$\\
  \end{enumerate}
\end{thm}

For the remainder of this article, we focus only on $g=0$ \GW\ invariants of a fixed smooth \DM\ stack $\sX$.  Thus, $\sX$ and $g$ will be suppressed in our notation.  In genus zero \GW\ theory, an important theorem is the WDVV equation, proven in the \DM\ stack context in \cite{A} and \cite{AGV2}.  To state it, consider the shorthand:  for any finite index set $I = \{1, \cdots, n\}$, and for any collection $\{\delta_i \in H^*(\ri_{\mu}(\sX))\}_{i \in I}$ we write $$\langle \gamma_1, \gamma_2, \delta_I, \ast \rangle_{\beta} \coloneq \langle \gamma_1, \gamma_2, \delta_{1}, \cdots, \delta_{n}, \ast \rangle_{\beta}.$$

\begin{thm}[\cite{A}, Theorem 6.2.1]
For all $\gamma_1, \gamma_2, \gamma_3 \in H^*(\ri_{\mu}(\sX))$, $\beta_3 \in N^+(X)$, and any finite collection $\{\delta_i \in H^*(\ri_{\mu}(\sX))\}_{i \in I}$ indexed by $I$, the following equation holds in $H^*(\ri_{\mu}(\sX))$:
\begin{alignat}{1}
&\sum_{A \sqcup B = I} \sum_{\beta_1 + \beta_2 = \beta_3} \langle \langle \gamma_1, \gamma_2, \delta_A, \ast \rangle_{\beta_1}, \gamma_3, \delta_B, \ast \rangle_{\beta_2}  \notag \\
= & \sum_{A \sqcup B = I} \sum_{\beta_1 + \beta_2 = \beta_3} \langle \langle \gamma_1, \gamma_3, \delta_A, \ast \rangle_{\beta_1}, \gamma_2, \delta_B, \ast \rangle_{\beta_2}. \notag
\end{alignat}
\end{thm}

In terms of a fixed additive basis $\alpha_0, \cdots, \alpha_p$ of $H^*(\ri_{\mu}(\sX))$ and corresponding matrix $\gi$ as in lemma (\ref{lemma}), the formulas can be expressed via \GW\ numbers:  for all $\gamma_1, \gamma_2, \gamma_3, \gamma_4$, the following equation holds in $R$:
\begin{alignat}{1} \label{WDVV}
& \sum_{A \sqcup B = I} \sum_{\beta_1 + \beta_2 = \beta_3} \sum_{i,j} \langle \gamma_1, \gamma_2, \delta_A, \alpha_i \rangle_{\beta_1}\ \gi\ \langle \alpha_j, \gamma_3, \delta_B, \gamma_4 \rangle_{\beta_2}  \notag \\
= & \sum_{A \sqcup B = I} \sum_{\beta_1 + \beta_2 = \beta_3} \sum_{i,j} \langle \gamma_1, \gamma_3, \delta_A, \alpha_i \rangle_{\beta_1}\ \gi\ \langle \alpha_j, \gamma_2, \delta_B, \gamma_4 \rangle_{\beta_2}.
\end{alignat}


\section{Quantum Cohomology}
Let $\sX$ be a smooth, proper \DM\ stack with projective coarse moduli scheme $X$.  The treatment here is standard; however, we shall desire to view the small quantum cohomology ring as a full deformation of $H^*(\ri_{\mu}(\sX), \bQ)$, so we shall make some convergence assumptions.
\begin{defn} \label{conv_crit}
$\sX$ satisfies the \emph{convergence criterion} if for each constant $C \in \bQ$, there are at most finitely many classes $\beta \in N^+(X)$ for which $\int_{\beta} c_1(T_{\sX}) < C$.
\end{defn}

As an example, any Fano scheme satisfies the convergence criterion (see \cite[Chapter~8]{CK}).

Let $\sX$ satisfy the convergence criterion, and let $\alpha_0, \cdots , \alpha_m$ be an additive basis for $H^*(\ri_{\mu}(\sX), \bQ)$ where $\alpha_0 = 1$ and $\alpha_1, \cdots, \alpha_p$ span $H^2(\sX, \bQ)$.  Let $\bQ[N^+(X)]$ be the monoid algebra of $N^+(X)$ over $\bQ$.  The algebra $\bQ[N^+(X)]$ has generators $\{\mathbf{q}^{\beta}\}_{\beta \in N^+(X)}$ and relations $\{\mathbf{q}^{\beta_1} \cdot \mathbf{q}^{\beta_2} - \mathbf{q}^{\beta_1 + \beta_2}\}_{\beta_1, \beta_2}$.

Consider the  $\bQ$-vector space given by
$$QH_{sm}^*(\sX) \coloneq H^*(\ri_{\mu},(\sX)) \otimes_{\bQ} \bQ[N^+(X)].$$

\begin{defn} \label{small_defn}
 The \emph{small quantum cohomology ring} is the vector space $QH_{sm}^*(\sX)$ with product given on the $\bQ [N^+(X)]$-basis by
$$\alpha_i \ast_{sm} \alpha_j = \sum_{\beta} \langle \alpha_i, \alpha_j, \ast \rangle_{\beta} \mathbf{q}^{\beta}.$$
\end{defn}

\begin{prop}
$QH^*(\sX)$ is a commutative, associative ring with identity.
\end{prop}
\begin{proof}
Associativity follows from the WDVV theorem (\ref{WDVV}).  Commutativity is clear from the definition.  That $\alpha_0$ is the identity follows from lemma (\ref{lemma}).
\end{proof}

One checks that this ring doesn't depend on the choice of basis of $H^*(\ri_{\mu}(\sX), \bQ)$.  Note that by the convergence criterion, the sum over $\beta$ in (\ref{small_defn}) is finite.  In some treatments, the convergence criterion is avoided by including formal variables corresponding to $\beta$ when defining the small quantum ring.

The ring $QH_{sm}^*(\sX)$ is regarded as a family of rings over $\bQ [N^+(X)]$.  We now describe the specializations.  For any point $\lambda \in \text{Spec } \bQ [N^+(X)]$, write $\bQ(\lambda)$ for the residue field at $\lambda$.

\begin{defn}
The \emph{quantum ring at $\lambda$}, denoted $QH^*_{\lambda}(\sX)$, is given by
$$QH^*_{\lambda}(\sX) \coloneq QH^*_{sm}(\sX) \otimes_{\bQ[N^+(X)]} \bQ(\lambda).$$
\end{defn}
Note that for each point $\lambda$, $QH^*_{\lambda}(\sX) \cong H^*(\ri_{\mu}(\sX), \bQ(\lambda))$ as $\bQ(\lambda)$-vector spaces, but not necessarily as rings.  Also, note that if $\bar{0}$ denotes the closed point corresponding to the ideal $\langle \mathbf{q}^{\beta} \rangle_{\beta \in N^+(X)}$, then $QH^*_{\bar{0}}(\sX)$ is the \emph{orbifold cohomology} or \emph{stringy cohomology ring} usually denoted as $H^*_{orb}(\sX)$ as defined in \cite{CR} and \cite{AGV}.  Furthermore, if $\sX \cong X$ is a smooth projective variety, then $QH_{\bar{0}}^*(X) \cong H^*(X, \bQ)$ is the usual cohomology ring.


\section{Reconstruction}

Suppose that $\sX$ is a smooth, proper \DM\ stack satisfying the convergence criterion.  The following is the main result of this article.

\begin{thm} \label{reconstruction}
  Suppose there exists $\lambda_0 \in \text{Spec }\bQ[N^+(X)]$ for which $QH^*_{\lambda_0}(\sX)$ is generated by untwisted divisor classes.  Then any genus zero \GW\ invariant on $\sX$ can be uniquely reconstructed from genus zero 3-point invariants.
\end{thm}
\begin{proof}
First we show that the condition of the hypothesis is an open condition on $\text{Spec } \bQ[N^+(X)]$.  Let $\alpha_1, \cdots, \alpha_p$ be a basis for the space $H^2(\sX, \bQ) \subset H^*(\ri_{\mu}(\sX), \bQ)$ of the untwisted divisor classes, and consider the $\bQ[N^+(X)]$-algebra homomorphism
$$\bQ[x_1, \cdots, x_p][N^+(X)] \xrightarrow{\phi} QH^*(\sX)$$
where the $x_i$ are indeterminants, and $\phi$ sends $x_i \mapsto \alpha_i$.  If $\phi|_{\bQ(\lambda)}$ denotes the induced morphism
$$\bQ[x_1, \cdots, x_p][N^+(X)] \otimes_{\bQ[N^+(X)]} \bQ(\lambda) \xrightarrow{\phi|_{\bQ(\lambda)}} QH^*_{\lambda}(\sX)$$
the condition is then equivalent to the surjectivity of $\phi|_{\bQ(\lambda_0)}$.  Regarding $\phi$ as a morphism of free $\sO_{\text{Spec }\bQ[N^+(X)]}$-modules, the locus where $\phi$ fails to be surjective is precisely $\text{Supp}(coker\ \phi)$, a closed subset since $coker\ \phi$ is a finitely generated module.

Now $\phi|_{\bQ(\lambda)}$ is surjective for all $\lambda$ in a open subset of $\text{Spec }\bQ[N^+(X)]$, so $\phi|_{\bQ(\eta)}$ is surjective where $\eta$ is the generic point of $\text{Spec }\bQ[N^+(X)]$.  Note that $N_1(X)$ is a finitely generated free abelian group, so $\text{Spec }\bQ[N^+(X)]$ is integral.  Thus $QH^*_{\eta}(\sX)$ is generated by untwisted divisor classes.

Next, consider computing the genus zero n-point \GW\ invariant
\begin{equation} \label{goal}
\langle \epsilon_1, \cdots, \epsilon_n \rangle_{\beta}
\end{equation}
with $\epsilon_i \in H^*(\ri_{\mu}(\sX), \bQ)$ and $\beta \in N^+(X)$.  Regard the $\epsilon_i$ as elements of $H^*(\ri_{\mu}(\sX), \bQ(\eta))$.  All the properties of \GW\ invariants from the previous sections apply to these classes as well.  By linearity we may assume $\epsilon_n = \alpha_1^{k_1} \cdot_{\eta} \ldots \cdot_{\eta} \alpha_p^{k_p}$.  Furthermore, we may assume $\sum k_i > 1$:  if $\sum k_i = 0$, (\ref{goal}) vanishes by the fundamental class axiom and if $\sum k_i = 1$, the untwisted divisor axiom reduces (\ref{goal}) to an (n - 1)-point invariant.  Thus for some $i$ with $k_i > 0$ write $D \coloneq \alpha_i$ and write $\epsilon_n = \epsilon_n' \cdot_{\eta} D$ for some $\epsilon_n'$.  Consider the WDWW equation (\ref{WDVV}) with $\gamma_1 = \epsilon_1$, $\gamma_2 = \epsilon_2$, $\gamma_3 = \epsilon_n'$, $\gamma_4 = D$, $\delta_I = \{\epsilon_3, \cdots, \epsilon_{n-1} \}$ and $\beta_3$ arbitrary.  The terms with $A = \emptyset$ or $B = \emptyset$ are exactly those involving 3-point invariants and n-point invariants.  Solving for these terms we can write (\ref{WDVV}) as
\begin{alignat}{1} \label{assoc}
\sum_{\beta_1 + \beta_2 = \beta_3} \langle \langle \epsilon_1, \epsilon_2, \ast \rangle_{\beta_1}, \epsilon_n', \delta_I, D \rangle_{\beta_2}
  &+  \sum_{\beta_1 + \beta_2 = \beta_3} \langle \langle \epsilon_1, \epsilon_2, \delta_I, \ast \rangle_{\beta_1}, \epsilon_n', D \rangle_{\beta_2} \notag \\
- \sum_{\beta_1 + \beta_2 = \beta_3} \langle \langle \epsilon_1, \epsilon_n', \ast \rangle_{\beta_1}, \epsilon_2, \delta_I, D \rangle_{\beta_2}
  &-  \sum_{\beta_1 + \beta_2 = \beta_3} \langle \langle \epsilon_1, \epsilon_n', \delta_I, \ast \rangle_{\beta_1}, \epsilon_2, D \rangle_{\beta_2} \\
  &= \text{terms with lower point invariants.} \notag
\end{alignat}
The untwisted divisor axiom applies to the first and third terms while the fourth term is calculated from an n-point invariant involving $\epsilon_n'$.  We shall induct on $n$ and $\sum k_i$.  Solving for the second term, the collection of the remaining terms on the right are determined by 3-point invariants via the inductive hypothesis.  Denoting these terms by $\Gamma(\beta_3)$, we obtain the following equation in $\bQ(\eta)$.

\begin{equation} \label{assoc2}
\sum_{\beta_1 + \beta_2 = \beta_3} \langle \epsilon_1, \cdots, \epsilon_{n-1}, \langle \epsilon_n', D, \ast \rangle_{\beta_2} \rangle_{\beta_1} = \Gamma(\beta_3)
\end{equation}
Now multiply (\ref{assoc2}) by $\mathbf{q}^{\beta_3}$ and sum over all $\beta_3$ to obtain

\begin{alignat}{2} \label{assoc3}
\sum_{\beta_3} \Gamma(\beta_3)\ \mathbf{q}^{\beta_3} \ \
    & &= &\sum_{\beta_3} \sum_{\beta_1 + \beta_2 = \beta_3} \langle \epsilon_1, \cdots, \epsilon_{n-1}, \langle \epsilon_n', D, \ast \rangle_{\beta_2}\ \mathbf{q}^{\beta_2} \rangle_{\beta_1}\ \mathbf{q}^{\beta_1} \notag \\
    & &= &\sum_{\beta_1} \langle \epsilon_1, \cdots, \epsilon_{n-1}, \sum_{\beta_2} \langle \epsilon_n', D, \ast \rangle_{\beta_2}\ \mathbf{q}^{\beta_2} \rangle_{\beta_1}\ \mathbf{q}^{\beta_1} \notag \\
    & &= &\sum_{\beta_1} \langle \epsilon_1, \cdots, \epsilon_{n-1}, \epsilon_n' \cdot_{\eta} D \rangle_{\beta_1}\ \mathbf{q}^{\beta_1} \notag \\
    & &= &\sum_{\beta_1} \langle \epsilon_1, \cdots, \epsilon_{n-1}, \epsilon_n \rangle_{\beta_1}\ \mathbf{q}^{\beta_1}. \notag
\end{alignat}
Thus we have the  following equality in $\bQ(\eta)$:
\begin{equation} \label{main}
\sum_{\beta_1} \langle \epsilon_1, \cdots, \epsilon_{n-1}, \epsilon_n \rangle_{\beta_1}\ \mathbf{q}^{\beta_1} = \sum_{\beta_3} \Gamma(\beta_3)\ \mathbf{q}^{\beta_3}.
\end{equation}
The right side of (\ref{main}) is determined inductively and is defined on some open set of $\text{Spec } \bQ[N^+(X)]$.  Hence it uniquely determines its coefficients on the left side:  the desired \GW\ invariants.
\end{proof}

\begin{rem}
Note that even when $QH^*_{\lambda_0}(\sX)$ is not generated by $H^2(\sX)$, we may restrict the \GW\ invariants considered by only pulling back classes from the subring generated by $H^2(\sX)$.  The theorem then reconstructs these invariants from the 3-point invariants.
\end{rem}


\section{Example: $\bP(1,b)$}

Let $\sX \cong \bP(1,b)$.  Recall the description of the rigidified cyclotomic inertia stack from Example \ref{mainexample}:
$$ \ri_{\mu}(\sX) = \sX \sqcup \bigsqcup_{1 \neq g \in \mu_b} B(\mu_b / \langle g \rangle).$$
Write $\mu_b \cong \bZ/b\bZ$, and label the component of $\ri_{\mu}(\sX)$ corresponding to $i \in \bZ/b\bZ$ by $B_i$:
$$ \ri_{\mu}(\sX) = \sX \sqcup \bigsqcup_{i = 1}^{b - 1} B_i.$$
Consider the additive basis for $H^*(\ri_{\mu}(\sX))$ given by $\alpha_0, \cdots, \alpha_{b-1}, x\ $ where $\alpha_0$ is the fundamental class of $\sX$, $\alpha_i$ is the fundamental class of $B_i$ for $1 \leq i \leq b-1$, and $x = c_1(\sO_{\sX}(1)) = [B(\mu_b)]$ is the fundamental class of the closed substack $B(\mu_b) \subset \sX$.  Finally, note that $N_1(X)$ is cyclic.  Let $\beta_0$ be a generator of $N^+(X)$, and let $\mathbf{q}$ denote the corresponding generator of $\bQ[N^+(X)]$.

The small quantum cohomology of $\sX$ is computed in \cite[Section~9]{AGV2} and is given by
\begin{equation} \label{QH_comp}
QH^*_{sm}(\sX) = \bQ[\mathbf{q}][\alpha_1]\ /\ \langle b \alpha_1^{b+1} - \mathbf{q} \rangle.
\end{equation}
Hence for the rest of the section we simply write $\alpha \coloneq \alpha_1$.

\begin{rem}
Note that in \cite{AGV2}, small quantum cohomology ring (denoted $QH^*(\sX)$) is defined as a formal deformation of $H^*(\ri_{\mu}(\sX))$ circumventing the need for any requirement on convergence of (\ref{small_defn}).  Since $\sX$ satisfies the convergence criterion, both rings are defined and $QH^*(\sX)$ is recovered from $QH^*_{sm}(\sX)$ by completing with respect to the ideal generated by $\mathbf{q}$.
\end{rem}

Now for any $\lambda \in \bQ$, we have
$$QH^*_{\lambda}(\sX) = \bQ[\alpha]\ /\ \langle b \alpha^{b+1} - \lambda \rangle$$
and one can check directly that for any $\lambda \neq 0$, $QH^*_{\lambda}(\sX)$ is generated by the untwisted divisor $x = \alpha^b$.  Theorem (\ref{reconstruction}) then applies.  Since $\bP(1,b)$ has Picard number $1$ and $\sX$ is \emph{Fano}, i.e.
$$\int_{\beta_0} c_1(T_{\bP(1,b)}) = \int_{\beta_0} (1 + \frac{1}{b})[pt] = 1 + \frac{1}{b} > 0,$$
the WDVV equation (\ref{assoc}) has an especially simple form.  For instance, consider computing an arbitrary degree 0, 4-point invariant on $\bP(1,b)$:
$$\langle \alpha^{k_1}, \alpha^{k_2}, \alpha^{k_3}, \alpha^{k_4} \rangle_0 \qquad \qquad k_1 + k_2 + k_3 + k_4 = 2b.$$
Proceed as in the proof of reconstruction.  Write $\alpha^{k_4} = \frac{b}{q}\alpha^{k_4 + 1} \ast \alpha^b$ so that $\langle \alpha^{k_4+1}, \alpha^b, \ast \rangle_1 = \frac{1}{b} \alpha^{k_4}$ (i.e. $D = \alpha^b$ and $\epsilon_n' = \alpha^{k_4 + 1}$).  Then (\ref{assoc}) becomes
\begin{alignat}{1} \label{assoc4}
\sum_{d_1 + d_2 = 1} \langle \langle \alpha^{k_1}, \alpha^{k_2}, \ast \rangle_{d_1}, &\alpha^{k_4+1}, \alpha^{k_3}, \alpha^b \rangle_{d_2} \notag \\
 &+ \sum_{d_1 + d_2 = 1} \langle \langle \alpha^{k_1}, \alpha^{k_2}, \alpha^{k_3}, \ast \rangle_{d_1}, \alpha^{k_4+1}, \alpha^b \rangle_{d_2} \notag \\
= \sum_{d_1 + d_2 = 1} \langle \langle \alpha^{k_1}, \alpha^{k_4+1}, \ast &\rangle_{d_1}, \alpha^{k_2}, \alpha^{k_3}, \alpha \rangle_{d_2} \\
 &-  \sum_{d_1 + d_2 = 1} \langle \langle \alpha^{k_1}, \alpha^{k_4+1}, \alpha^{k_3}, \ast \rangle_{d_1}, \alpha^{k_2}, \alpha^b \rangle_{d_2} \notag
\end{alignat}
By linearity, we may rewrite the second and fourth terms above.  For instance in the second term
$$\langle \langle \alpha^{k_1}, \alpha^{k_2}, \alpha^{k_3}, \ast \rangle_{d_1}, \alpha^{k_4+1}, \alpha^b \rangle_{d_2} = \langle \alpha^{k_1}, \alpha^{k_2}, \alpha^{k_3}, \langle \alpha^{k_4+1}, \alpha^b , \ast \rangle_{d_2} \rangle_{d_1}.$$
This allows us to simplify (\ref{assoc4}) using explicit knowledge of the small quantum product.  Since (\ref{QH_comp}) yields
\begin{equation} \label{explicit_comp}
\langle \alpha^{k_1}, \alpha^{k_2}, \ast \rangle_d =
  \begin{cases}
  \alpha^{k_1 + k_2} & \text{if $d = 0$ and $k_1 + k_2 \leq b$}\\
  \frac{b}{q}\alpha^{k_1 + k_2 - b} & \text{if $d = 1$ and $k_1 + k_2 > b$}\\
  0 & \text{otherwise}
  \end{cases}
\end{equation}
the first term in (\ref{assoc4}) admits two cases.  First if $k_1 + k_2 > b$, then by (\ref{explicit_comp}) the only contribution occurs when $d_1 = 1$.  This term then becomes
$$\frac{1}{b} \langle \alpha^{k_1 + k_2 - b}, \alpha^{k_4+1}, \alpha^{k_3}, \alpha^b \rangle_0$$
which must vanish by the Degree Axiom.  Second if $k_1 + k_2 \leq b$, then the only contribution occurs when $d_1 = 0$.  This term then becomes
$$\langle \alpha^{k_1 + k_2}, \alpha^{k_4+1}, \alpha^{k_3}, \alpha^b \rangle_1$$
which equals
$$\frac{1}{b} \int_{\beta_0} \alpha^b = \frac{1}{b^3}$$
by the Divisor Axiom and (\ref{explicit_comp}).  Thus we define the function $e(l)$ to be $0$ if $l > b$ and $1$ otherwise.  The first term in (\ref{assoc4}) then simplifies to $\frac{e(k_1 + k_2)}{b^3}$.  The remaining terms simplify similarly and (\ref{assoc4}) becomes
$$\frac{e(k_1 + k_2)}{b^3} + \frac{1}{b}\langle \alpha^{k_1}, \alpha^{k_2}, \alpha^{k_3}, \alpha^{k_4} \rangle_0  = \frac{e(k_1 + k_4 + 1)}{b^3} +  \frac{1}{b} \langle \alpha^{k_1}, \alpha^{k_4 + 1}, \alpha^{k_3}, \alpha^{k_2 -1} \rangle_0.$$
It is clear that we can iterate this process increasing the exponent of $\alpha^{k_4+1}$ and decreasing the exponent of $\alpha^{k_2-1}$.  Assuming we have chosen $k_4$ as the largest exponent and $k_2$ as the next largest, we may iterate exactly $b-k_4$ times.  Then $\alpha^{k_4 + (b-k_4)} = \alpha^b$ is the divisor class and we obtain the formula
\begin{equation} \label{formula_1}
\langle \alpha^{k_1}, \alpha^{k_2}, \alpha^{k_3}, \alpha^{k_4} \rangle_0 = \frac{(b-(k_1+k_4))e(k_1+k_4) - k_3}{b^2}
\end{equation}
Making different choices to apply the above algorithm gives other compuations similar to (\ref{formula_1}).  By averaging these computations we obtain a formula independent of the various choices made:
\begin{equation} \label{formula_2}
\langle \alpha^{k_1}, \alpha^{k_2}, \alpha^{k_3}, \alpha^{k_4} \rangle_0 = \frac{\sum_{i \neq j} (b-(k_i+k_j))e(k_i+k_j) - b}{2b^2}
\end{equation}

Similarly we obtain
\begin{equation} \label{formula_3}
\langle \alpha^{k_1}, \alpha^{k_2}, \alpha^{k_3}, \alpha^{k_4} \rangle_1 = \frac{1}{b^3} \qquad \qquad k_1 + k_2 + k_3 + k_4 = 3b + 1.
\end{equation}

Closed formulae for higher point invariants are more complicated, however the reconstruction theorem still gives an explicit algorithm for computing an arbitrary genus zero \GW\ invariant on $\bP(1,b)$.  This algorithm was implemented in Maple and the appendix contains complete calculations for small $b$.

\begin{rem}
The genus zero \GW\ invariants of more general weighted projective spaces are investigated by Mann \cite{Mann}, and more recently by Coates, Corti, Lee and Tseng \cite{CCLT}.
\end{rem}


\appendix
\section{$\bP(1,b)$ Computations}
Let $\sX \cong \bP(1,b)$ and let $\alpha_0, \alpha_1, \cdots, \alpha_{b-1}, x$ be the additive basis of $H^*(\ri_{\mu}(\sX))$ as in Section 4.  Consider the shorthand $N_d(k_1, \cdots, k_{b-1}) \coloneq \langle\ \underbrace{\alpha_1, \cdots, \alpha_1}_{k_1}, \cdots, \underbrace{\alpha_{b-1}, \cdots, \alpha_{b-1}}_{k_{b-1}}\ \rangle_{d \beta_0}$.  Omitting invariants with an untwisted divisor, the following computations tabulate all non-zero $n > 3$ invariants for $\bP(1,b), b = 1, \cdots, 6$.

\begin{table}[h]
  \begin{flushleft}
    \renewcommand{\arraystretch}{\TABLESIZE}
    \begin{tabular}{lll}
        $\sX = \bP(1,2)$ &\hspace{6.0cm} & $\sX = \bP(1,3)$\\
        \cline{1-1} \cline{3-3}
        $N_0(4) = -1/4$ &\hspace{6.0cm} &  $N_0(2,2) = -1/9$\\
                        &\hspace{6.0cm} &  $N_0(1,4) = 1/27$\\
                        &\hspace{6.0cm} &  $N_0(0,6) = -1/27$\\
    \end{tabular}
  \end{flushleft}
\end{table}
\begin{table}[h]
  \begin{flushleft}
    \renewcommand{\arraystretch}{\TABLESIZE}
    \begin{tabular}{l|l}
        \multicolumn{2}{l}{$\sX = \bP(1,4)$} \\
        \hline
                     $N_0(1,2,1) = -1/16$
                   & $N_0(0,4,0) = -1/8$ \\
                     $N_0(2,0,2) = -1/16$
                   & $N_0(1,1,3) = 1/64$ \\
                     $N_0(0,3,2) = 1/32$
                   & $N_0(0,2,4) = -1/64$ \\
                     $N_0(1,0,5) = -1/256$
                   & $N_0(0,1,6) = 5/512$ \\
                     $N_0(0,0,8) = -5/512$ \\
    \end{tabular}
  \end{flushleft}
\end{table}
\begin{table}[h]
  \begin{center}
    \renewcommand{\arraystretch}{\TABLESIZE}
    \begin{tabular}{l|l|l}
        \multicolumn{3}{l}{$\sX = \bP(1,5)$} \\
        \hline
                    $N_0(1,1,1,1) = -1/25$
                  & $N_0(1,0,2,2) = 1/125$
                  & $N_0(0,2,1,2) = 2/125$ \\
                    $N_0(1,1,0,3) = 1/125$
                  & $N_0(0,1,2,3) = -1/125$
                  & $N_0(0,1,3,1) = 3/125$ \\
                    $N_0(1,0,3,0) = -1/25$
                  & $N_0(1,0,1,4) = -1/625$
                  & $N_0(0,2,2,0) = -2/25$ \\
                    $N_0(2,0,0,2) = -1/25$
                  & $N_0(0,3,0,1) = -1/25$
                  & $N_0(0,2,0,4) = -4/625$ \\
                    $N_0(0,1,1,5) = 11/3125$
                  & $N_0(0,0,3,4) = 21/3125$
                  & $N_0(0,0,4,2) = -9/625$ \\
                    $N_0(0,0,2,6) = -13/3125$
                  & $N_0(0,1,0,7) = -1/625$
                  & $N_0(1,0,0,6) = 1/3125$ \\
                    $N_0(0,0,1,8) = 49/15625$
                  & $N_0(0,0,5,0) = 1/25$
                  & $N_0(0,0,0,10) = -49/15625$ \\
                    $N_1(0,0,0,4) = 1/125$ \\
    \end{tabular}
  \end{center}
\end{table}
\begin{table}[!h]
  \begin{center}
    \renewcommand{\arraystretch}{\TABLESIZE}
    \begin{tabular}{l|l|l}
        \multicolumn{3}{l}{$\sX = \bP(1,6)$} \\
        \hline
                    $N_0(0,0,3,0,3) = -1/216$
                  & $N_0(0,0,1,2,5) = -91/46656$
                  & $N_0(0,1,0,1,6) = -13/23328$ \\
                   $N_0(0,1,1,0,5) = 11/7776$
                  & $N_0(0,0,2,1,4) = 11/3888$
                  & $N_0(0,1,0,2,4) = 1/648$ \\
                   $N_0(0,2,0,0,4) = -1/324$
                  & $N_0(0,1,0,3,2) = -1/216$
                  & $N_0(0,0,3,1,1) = 1/54$ \\
                   $N_0(1,0,0,1,5) = 1/7776$
                  & $N_0(0,0,2,2,2) = -5/648$
                  &$N_0(0,2,0,2,0) = -1/18$ \\
                   $N_0(1,0,1,0,4) = -1/1296$
                  & $N_0(1,0,0,2,3) = -1/1296$
                  &$N_0(1,0,0,3,1) = 1/216$ \\
                   $N_0(0,1,2,0,2) = 1/108$
                  & $N_0(0,1,2,1,0) = -1/18$
                  & $N_0(0,2,0,1,2) = 1/108$ \\
                   $N_0(0,1,1,1,3) = -5/1296$
                  & $N_0(0,2,1,0,1) = -1/36$
                  & $N_0(1,1,0,0,3) = 1/216$ \\
                   $N_0(0,1,1,2,1) = 1/72$
                  & $N_0(2,0,0,0,2) = -1/36$
                  & $N_0(1,0,2,0,1) = -1/36$ \\
                   $N_0(1,0,1,2,0) = -1/36$
                  & $N_0(1,0,1,1,2) = 1/216$
                  & $N_0(1,1,0,1,1) = -1/36$ \\
                   $N_0(0,0,4,0,0) = -1/12$
                  & $N_0(0,0,2,0,6) = -11/7776$
                  &$N_0(0,0,2,3,0) = 1/36$\\
                   $N_0(0,1,0,0,8) = 7/34992$
                  & $N_0(0,1,0,4,0) = 1/54$
                  &$N_0(1,0,0,0,7) = -1/46656$\\
                   $N_0(0,0,1,1,7) = 301/279936$
                  & $N_0(0,0,1,3,3) = 11/2592$
                  &$N_0(0,0,1,4,1) = -5/432$\\
                   $N_0(0,0,1,0,9) = -119/186624$
                  & $N_0(0,0,0,3,6) = 91/46656$
                  &$N_0(0,0,0,4,4) = -13/3888$\\
                   $N_0(0,0,0,6,0) = -1/54$
                  & $N_0(0,0,0,5,2) = 1/144$
                  &$N_0(0,0,0,2,8) = -287/209952$\\
                   $N_0(0,0,0,1,10) = 5663/5038848$
                  & $N_0(0,0,0,0,12) = -5663/5038848$
                  &$N_1(0,0,0,1,3) = 1/216$\\
                   $N_1(0,0,0,0,5) = 1/1296$
    \end{tabular}
  \end{center}
\end{table}



\begin{thebibliography}{99}
\bibitem{A}
    D. Abramovich,
    \emph{Lectures on Gromov-Witten invariants of orbifolds},
    \newline arXiv:math.AG/0512372.

\bibitem{AGV}
    D. Abramovich, T. Graber, and A. Vistoli,
    \emph{Algebraic orbifold quantum products},
    Orbifolds in mathematics and physics (Madison, WI, 2001),
    Contemp. Math.,
    vol. 310,
    Amer. Math. Soc.,
    Providence, RI,
    2002,
    pp. 1-24.

\bibitem{AGV2}
    D. Abramovich, T. Graber, and A. Vistoli,
    \emph{Gromov--Witten theory of Deligne--Mumford stacks},
    arXiv:math.AG/0603151.

\bibitem{AV1}
    D. Abramovich and A. Vistoli,
    \emph{Compactifying the space of stable maps},
    J. Amer. Math. Soc. \textbf{15} (2002),
    no. 1, 27-75.

\bibitem{BF}
    K. Behrend and B. Fantechi,
    \emph{The intrinsic normal cone},
    Invent. Math. \textbf{128} (1997),
    no. 1, 45-88.

\bibitem{Be}
    K. Behrend and Yu. Manin,
    \emph{Stacks of stable maps and Gromov-Witten invariants},
    Duke Math. J. \textbf{85} (1996),
    no. 1, 1-60.

\bibitem{CR}
    W. Chen and Y. Ruan,
    \emph{A new cohomology theory of orbifolds},
    Comm. Math. Phys. \textbf{248} (2004),
    no. 1, 1-31.

\bibitem{CCLT}
    T. Coates, A. Corti, Y.-P. Lee and H.-H. Tseng,
    \emph{The quantum orbifold cohomology of weighted projective space},
    arXiv:math.AG/0608481.

\bibitem{CK}
    D. Cox and S. Katz,
    \emph{Mirror symmetry and algebraic geometry},
    Mathematical Surveys and Monographs,
    vol. 68,
    American Mathematical Society,
    Providence, RI,
    1999.

\bibitem{DVV2}
    R. Dijkgraaf, H. Verlinde and E. Verlinde,
    \emph{Notes on topological string theory and $2$D quantum
              gravity},
    String theory and quantum gravity (Trieste, 1990),
    World Sci. Publ., River Edge, NJ,
    1991, pp. 91-156.
    
\bibitem{DVV1}
    R. Dijkgraaf, H. Verlinde and E. Verlinde,
    \emph{Topological strings in $d<1$},
    Nuclear Phys. B \textbf{352} (1991),
    no. 1, 59-86.

\bibitem{Ko}
    M. Kontsevich,
    \emph{Enumeration of rational curves via torus actions},
    The moduli space of curves (Texel Island, 1994),
    Progr. Math.,
    vol. 129,
    Birkh\"auser Boston, Boston, MA,
    1995, pp. 335-368.

\bibitem{KM}
    M. Kontsevich and Yu. Manin,
    \emph{Gromov-Witten classes, quantum cohomology, and enumerative
              geometry},
    Mirror symmetry, II,
    AMS/IP Stud. Adv. Math.,
    vol. 1,
    Amer. Math. Soc.,
    Providence, RI,
    1997, pp. 607-653.

\bibitem{Mann}
    E. Mann,
    \emph{Cohomologie quantique orbifolde des espaces projectifs a
        poids},
    \newline arXiv:math.AG/0510331.

\bibitem{O}
    M. Olsson,
    \emph{On log twisted curves},
    \newline \textbf{http://www.ma.utexas.edu/users/molsson/Logcurves.pdf}.

\bibitem{O-Hom}
    M. Olsson,
    \emph{$\underline {\rm Hom}$-stacks and restriction of scalars},
    Duke Math. J. \textbf{134} (2006),
    no. 1, 139-164.

\bibitem{Witten}
    E. Witten,
    \emph{On the structure of the topological phase of two-dimensional
              gravity},
    Nuclear Phys. B, \textbf{340} (1990),
    no. 2-3, 281-332.
\end{thebibliography}
\end{document}